\documentclass[12pt,thmsa]{article}
\usepackage{sw20bams}

\input tcilatex
\begin{document}

\author{Steven R. Finch}
\title{Another Look at AR(1)}
\date{December 29, 2007}
\maketitle

\begin{abstract}
Given a stationary first-order autoregressive process $X_t$ (with lag-one
correlation $\rho $ satisfying $|\rho |<1$), we examine the Central Limit
Theorem for $\frac 1n\ln |X_1\cdots X_n|$ and compute variances to high
precision. Given a nonstationary process $X_t$ (with $|\rho |>1$), we
examine instead $\frac 1n\ln |X_n|$ and study the distribution of $\ln
|X_n|-n\,\ln |\rho |$.
\end{abstract}

\footnotetext{
Copyright \copyright\ 2007 by Steven R. Finch. All rights reserved.}This
research began as an effort to better understand Viswanath's random integer
recurrence \cite{Vi}: 
\[
\begin{array}{ccccc}
X_t=X_{t-1}\pm X_{t-2}, &  & X_0=1, &  & X_1=1,
\end{array}
\]
\[
\begin{array}{ccc}
\frac 1n\ln |X_n|\rightarrow \ln (1.13198824...) &  & \text{almost surely as 
}n\rightarrow \infty
\end{array}
\]
and one of Wright \&\ Trefethen's real recurrences \cite{WT}: 
\[
\begin{array}{ccccc}
X_t=X_{t-1}+\varepsilon _tX_{t-2}, &  & X_0=1, &  & X_1=1,
\end{array}
\]
\[
\begin{array}{ccc}
\frac 1n\ln |X_n|\rightarrow \ln (1.057473553704...) &  & \text{almost
surely as }n\rightarrow \infty
\end{array}
\]
where $\varepsilon _t$ is $N(0,1)$ white noise. What are the corresponding
asymptotic results for certain well-known recurrences in standard time
series analysis?

When $|\rho |>1$, the nonstationary first-order autoregressive process 
\[
\begin{array}{ccc}
X_t=\rho \,X_{t-1}+\sqrt{\rho ^2-1}\,\varepsilon _t, &  & X_0=0
\end{array}
\]
is readily shown to satisfy 
\[
\begin{array}{ccc}
\frac 1n\ln |X_n|\rightarrow \ln |\rho | &  & \text{almost surely as }%
n\rightarrow \infty .
\end{array}
\]
The quantity $\ln |\rho |$ is called the \textit{Lyapunov exponent} of the
system \cite{CT}. More precisely, 
\[
\begin{array}{ccc}
\mu _n=\limfunc{E}(\ln |X_n|)=\frac 12\left( \ln (\rho ^{2n}-1)-\ln
(2)-\gamma \right) , &  & \sigma ^2=\limfunc{Var}(\ln |X_n|)=\frac{\pi ^2}8
\end{array}
\]
where $\gamma $ denotes the Euler-Mascheroni constant \cite{Fi1}. We wish to
ascertain the distribution of the errors $\left( \ln |X_n|-\mu _n\right)
/\sigma $, which do not appear to be $N(0,1)$.

When $|\rho |<1$, the stationary first-order autoregressive process 
\[
X_t=\rho \,X_{t-1}+\sqrt{1-\rho ^2}\,\varepsilon _t 
\]
gives rise to a different question. Here we have 
\[
\begin{array}{ccc}
\mu =\limfunc{E}(\ln |X_t|)=\frac 12\left( -\ln (2)-\gamma \right) , &  & 
\sigma ^2=\limfunc{Var}(\ln |X_t|)=\frac{\pi ^2}8
\end{array}
\]
in contrast to before. The Central Limit Theorem gives \cite{S1, S2} 
\[
\begin{array}{ccc}
\sqrt{n}\dfrac{\frac 1n\dsum\limits_{t=1}^n\ln |X_t|-\mu }{\xi _\rho }%
\rightarrow N(0,1) &  & \text{as }n\rightarrow \infty
\end{array}
\]
for some constant $\xi _\rho >0$; clearly $\xi _0=\sigma $. What is the
numerical value of $\xi _\rho $ as a function of $\rho \neq 0$? This is our
first question to be addressed.

\subsection{Stationary Case}

Let $f(x)$ denote the $N(0,1)$ density function and $f^{(j)}(x)$ denote its $%
j^{\text{th}}$ derivative. Since $\limfunc{Cov}\left( X_1,X_{\ell +1}\right)
=\rho ^\ell $ for integer lag $\ell \geq 1$, it follows that \cite{Ne} 
\begin{eqnarray*}
\limfunc{E}\left( \ln |X_1|\cdot \ln |X_{\ell +1}|\right)
&=&\dsum\limits_{j=0}^\infty \left| \dint\limits_{-\infty }^\infty \ln
|x|\,f^{(j)}(x)\,dx\right| ^2\frac{\rho ^{j\,\ell }}{j!} \\
&=&\mu ^2+\dsum\limits_{k=1}^\infty \nu _{2k}^2\frac{\rho ^{2k\,\ell }}{(2k)!%
}
\end{eqnarray*}
where 
\[
\nu _{2k}=(-1)^{k-1}\dint\limits_{-\infty }^\infty \ln
|x|\,f^{(2k)}(x)\,dx=2^{k-1}(k-1)! 
\]
Hence 
\begin{eqnarray*}
\limfunc{Var}\left( \frac 1{\sqrt{n}}\dsum\limits_{t=1}^n\ln |X_t|\right)
&=&\frac 1n\dsum\limits_{t=1}^n\dsum\limits_{s=1}^n\limfunc{Cov}\left( \ln
|X_t|,\ln |X_s|\right) \\
&=&\sigma ^2+\frac 2n\dsum\limits_{\ell =1}^{n-1}(n-\ell )\left( \limfunc{E}%
\left( \ln |X_1|\cdot \ln |X_{\ell +1}|\right) -\mu ^2\right) \\
&=&\sigma ^2+\frac 2n\dsum\limits_{k=1}^\infty \frac{\nu _{2k}^2}{(2k)!}%
\dsum\limits_{\ell =1}^{n-1}(n-\ell )\rho ^{2k\,\ell } \\
&=&\sigma ^2+\frac 2n\dsum\limits_{k=1}^\infty \frac{\nu _{2k}^2}{(2k)!}%
\left( \frac n{1-\rho ^{2k}}-\frac{1-\rho ^{2k\,n}}{(1-\rho ^{2k})^2}\right)
\rho ^{2k}
\end{eqnarray*}
and therefore 
\begin{eqnarray*}
\xi _\rho ^2 &=&\lim_{n\rightarrow \infty }\limfunc{Var}\left( \frac 1{\sqrt{%
n}}\dsum\limits_{t=1}^n\ln |X_t|\right) =\sigma
^2+2\dsum\limits_{k=1}^\infty \frac{\nu _{2k}^2}{(2k)!}\frac{\rho ^{2k}}{%
1-\rho ^{2k}} \\
&=&\frac{\pi ^2}8+2\left( \frac 12\frac{\rho ^2}{1-\rho ^2}+\frac 16\frac{%
\rho ^4}{1-\rho ^4}+\frac 4{45}\frac{\rho ^6}{1-\rho ^6}+\frac 2{35}\frac{%
\rho ^8}{1-\rho ^8}\right. \\
&&\left. +\frac{64}{1575}\frac{\rho ^{10}}{1-\rho ^{10}}+\frac{64}{2079}%
\frac{\rho ^{12}}{1-\rho ^{12}}+\frac{512}{21021}\frac{\rho ^{14}}{1-\rho
^{14}}+\frac{128}{6435}\frac{\rho ^{16}}{1-\rho ^{16}}\right. \\
&&\left. +\frac{16384}{984555}\frac{\rho ^{18}}{1-\rho ^{18}}+\frac{16384}{%
1154725}\frac{\rho ^{20}}{1-\rho ^{20}}+\frac{131072}{10669659}\frac{\rho
^{22}}{1-\rho ^{22}}+\cdots \right)
\end{eqnarray*}
via computer algebra. This is an example of what is called a \textit{Lambert
series }\cite{BB}. With suitably many terms, we calculate 
\[
\xi _{0.1}=1.11527354305263680232..., 
\]
\[
\xi _{0.3}=1.15562165351986837602..., 
\]
\[
\xi _{0.5}=1.26199222423122947973..., 
\]
\[
\xi _{0.7}=1.52783735828651737636..., 
\]
\[
\xi _{0.9}=2.55564072887132125752... 
\]
to 20 decimal places.

As a corollary, if $Y_t$ is an Ornstein-Uhlenbeck process (Gauss-Markov
process) satisfying 
\[
\begin{array}{ccc}
dY_t=-\theta \,Y_t\,dt+\sqrt{2\theta }\,dW_t, &  & 0\leq t\leq T
\end{array}
\]
where $\theta >0$ and $W_t$ is Brownian motion with unit variance, then \cite
{S1} 
\[
\begin{array}{ccc}
\sqrt{T}\dfrac{\frac 1T\dint\limits_0^T\ln \left| Y_t\right| dt-\mu }{\eta
_\theta }\rightarrow N(0,1) &  & \text{as }T\rightarrow \infty
\end{array}
\]
for some constant $\eta _\theta >0$. A formula for $\eta _\theta $ is proved
as follows \cite{Ne}:\ 
\begin{eqnarray*}
\limfunc{E}\left( \ln |Y_0|\cdot \ln |Y_\ell |\right)
&=&\dsum\limits_{j=0}^\infty \left| \dint\limits_{-\infty }^\infty \ln
|y|\,f^{(j)}(y)\,dy\right| ^2\frac{e^{-j\,\theta \,\ell }}{j!} \\
\ &=&\mu ^2+\dsum\limits_{k=1}^\infty \nu _{2k}^2\frac{e^{-2k\,\theta \,\ell
}}{(2k)!}
\end{eqnarray*}
because $\limfunc{Cov}\left( Y_0,Y_\ell \right) =e^{-\theta \,\ell }$ for
real lag $\ell >0$; hence 
\begin{eqnarray*}
\limfunc{Var}\left( \frac 1{\sqrt{T}}\dint\limits_0^T\ln |Y_t|\,dt\right)
&=&\frac 1T\dint\limits_0^T\dint\limits_0^T\limfunc{Cov}\left( \ln |Y_t|,\ln
|Y_s|\right) \,ds\,dt \\
\ &=&\frac 2T\dint\limits_0^T\dint\limits_0^t\limfunc{Cov}\left( \ln
|Y_t|,\ln |Y_{t-\ell }|\right) \,d\ell \,dt
\end{eqnarray*}
upon setting $\ell =t-s$, $d\ell =-ds$ for fixed $t$; hence 
\[
\limfunc{Var}\left( \frac 1{\sqrt{T}}\dint\limits_0^T\ln |Y_t|\,dt\right)
=\frac 2T\dint\limits_0^T\dint\limits_\ell ^T\limfunc{Cov}\left( \ln
|Y_{t-\ell }|,\ln |Y_t|\right) \,dt\,d\ell 
\]
upon reversing the order of integration; hence 
\begin{eqnarray*}
\limfunc{Var}\left( \frac 1{\sqrt{T}}\dint\limits_0^T\ln |Y_t|\,dt\right)
&=&\frac 2T\dint\limits_0^T\dint\limits_\ell ^T\limfunc{Cov}\left( \ln
|Y_0|,\ln |Y_\ell |\right) \,dt\,d\ell \\
\ &=&\frac 2T\dint\limits_0^T\left( T-\ell \right) \left( \limfunc{E}\left(
\ln |Y_0|\cdot \ln |Y_\ell |\right) -\mu ^2\right) \,d\ell \\
\ &=&\frac 2T\dsum\limits_{k=1}^\infty \frac{\nu _{2k}^2}{(2k)!}%
\dint\limits_0^T\left( T-\ell \right) e^{-2k\,\theta \,\ell }\,d\ell \\
\ &=&\frac 2T\dsum\limits_{k=1}^\infty \frac{\nu _{2k}^2}{(2k)!}\left( \frac
T{2k\,\theta }+\frac{e^{-2k\,\theta \,T}-1}{(2k\,\theta )^2}\right)
\end{eqnarray*}
and therefore 
\begin{eqnarray*}
\eta _\theta ^2 &=&\lim_{T\rightarrow \infty }\limfunc{Var}\left( \frac 1{%
\sqrt{T}}\dint\limits_0^T\ln |Y_t|\,dt\right) =\frac 1\theta
\dsum\limits_{k=1}^\infty \frac{\nu _{2k}^2}{(2k)!}\frac 1k \\
\ &=&\frac 1\theta \dsum\limits_{k=1}^\infty \frac{2^{2k-2}(k-1)!^2}{(2k)!\,k%
}=\frac 1\theta \dsum\limits_{n=0}^\infty \frac{2^{2n-1}}{\tbinom{2n}%
n(n+1)^2(2n+1)} \\
\ &=&\frac 1\theta \left( \frac 14\pi ^2\ln (2)-\frac 78\zeta (3)\right)
=\frac 1\theta (0.81146307722510340753...)^2
\end{eqnarray*}
where $\zeta (s)$ denotes the Riemann zeta function \cite{Fi2}. Many
analogous central binomial sums appear in \cite{Sh}.

\subsection{Nonstationary Case}

Since $X_n\sim N(0,\rho ^{2n}-1)$, we deduce that 
\begin{eqnarray*}
\limfunc{P}\left( \frac{\ln |X_n|-\mu _n}\sigma \leq x\right) &=&\limfunc{P}%
\left( |X_n|\leq e^{\sigma \,x+\mu _n}\right) \\
&=&\sqrt{\frac 2{\pi (\rho ^{2n}-1)}}\dint\limits_0^{e^{\sigma \,x+\mu
_n}}\exp \left( -\frac{y^2}{2(\rho ^{2n}-1)}\right) dy;
\end{eqnarray*}
thus 
\begin{eqnarray*}
\frac d{dx}\limfunc{P}\left( \frac{\ln |X_n|-\mu _n}\sigma \leq x\right) &=&%
\sqrt{\frac 2{\pi (\rho ^{2n}-1)}}\exp \left( -\frac{e^{2(\sigma \,x+\mu _n)}%
}{2(\rho ^{2n}-1)}\right) e^{\sigma \,x+\mu _n}\sigma \\
&=&\frac 12\sqrt{\frac \pi {\rho ^{2n}-1}}\exp \left( -\frac{e^{2(\sigma
\,x+\mu _n)}}{2(\rho ^{2n}-1)}+\sigma \,x+\mu _n\right) \\
&=&\frac 12\sqrt{\frac \pi 2}\exp \left( -\frac 14e^{2\sigma \,x-\gamma
}+\sigma \,x-\frac \gamma 2\right) \\
&=&\frac 12\sqrt{\frac \pi 2}\exp \left( -\frac 14e^z+\frac 12z\right)
\end{eqnarray*}
where $z=\pi x/\sqrt{2}-\gamma $. Clearly $\left( \ln |X_n|-\mu _n\right)
/\sigma $ possesses a doubly exponential density function (called a Gumbel
density or Fisher-Tippett Type I extreme values density \cite{Re}) with mean 
$=0$, variance $=1$, 
\[
\text{skewness}=-\dfrac{28\sqrt{2}}{\pi ^3}\zeta
(3)=-1.53514159072290597506... 
\]
and kurtosis $=7-3=4$. Negativity of the third moment above confirms that
the distribution is skewed to the left. Closed-form expressions for the
quartiles do not exist: 
\[
25^{\text{th}}\text{ }\%\text{-tile}=-0.45782337329420373497..., 
\]
\[
\text{median}=50^{\text{th}}\text{ }\%\text{-tile}%
=0.21732071404060381038..., 
\]
\[
75^{\text{th}}\text{ }\%\text{-tile}=0.69796763838144042777... 
\]
but the maximum point of the density is easily found: 
\[
\text{mode}=\frac{\sqrt{2}(\ln (2)+\gamma )}\pi =0.57186419860436852975.... 
\]
It is pleasing that, upon subtracting the ``trend'' from an AR(1) process,
such a nice residual distribution emerges (independent of both $\rho $ and $%
n $).

As a corollary, if $Y_t$ satisfies 
\[
\begin{array}{ccccc}
dY_t=-\theta \,Y_t\,dt+\sqrt{-2\theta }\,dW_t, &  & Y_0=0, &  & 0\leq t\leq T
\end{array}
\]
where $\theta <0$, then 
\[
\begin{array}{ccc}
\frac 1T\ln |Y_T|\rightarrow -\theta &  & \text{almost surely as }%
T\rightarrow \infty
\end{array}
\]
and the density of $\left( \ln |Y_T|+T\,\theta \right) /\sigma $ approaches
the same doubly exponential function as before. The proof is immediate.

More generally, consider the nonstationary AR($m$) process 
\[
X_t=a_1X_{t-1}+a_2X_{t-2}+\cdots +a_mX_{t-m}+b\,\varepsilon _t, 
\]
\[
X_0=X_{-1}=\cdots =X_{2-m}=X_{1-m}=0. 
\]
Let $A$ denote the $m\times m$ matrix with $\left( a_1,a_2,\ldots
,a_m\right) $ in the top row, $1$s on the subdiagonal and $0$s elsewhere.
Order the complex eigenvalues $\lambda _1$, $\lambda _2$, $\ldots $, $%
\lambda _m$ of $A$ so that $\lambda _1$ has maximum modulus. When $|\lambda
_1|>1$, AR($m$) is shown to satisfy \cite{ME} 
\[
\begin{array}{ccc}
\frac 1n\ln |X_n|\rightarrow \ln |\lambda _1| &  & \text{almost surely as }%
n\rightarrow \infty .
\end{array}
\]
This occurs, for $m=2$, if and only if $|a_2|>1$ or $|a_1|>1-a_2$. An
evaluation of the residual distribution remains open.

\subsection{Variations}

Surely the results given in this paper are not new! A\ careful literature
search was unsuccessful. An example in \cite{CT} inspires us to look at the
stationary case with $\varepsilon _t$ assumed to be $U(-\sqrt{3},\sqrt{3})$
white noise. Obviously $\limfunc{E}(X_t)=0$ and $\limfunc{Var}(X_t)=1$. When 
$\rho =0$, it follows that 
\[
\begin{array}{ccc}
\limfunc{E}(\ln |X_t|)=\tfrac 12\ln (3)-1\approx \ln (0.637), &  & \limfunc{%
Var}(\ln |X_t|)=1
\end{array}
\]
because each $X_t$ is uniformly distributed.\footnote[1]{%
The numerical estimate $\limfunc{E}(\ln |X_t|)\approx \ln (0.2)$ in \cite{CT}
is evidently a mistake.} When $\rho \neq 0$, this fact no longer holds and
hence the relevant Central Limit Theorem parameters are not apparent.

We conclude with a recurrence that somewhat resembles Viswanath's: 
\[
\begin{array}{ccc}
X_t=\rho \,X_{t-1}\pm \sqrt{\rho ^2-1}, &  & X_0=0
\end{array}
\]
where $|\rho |>1$ and plus/minus signs are equiprobable. While $\limfunc{E}%
(X_n)=0$ and $\limfunc{Var}(X_n)=\rho ^{2n}-1$ as in the Gaussian
nonstationary case, it seems difficult to find $\limfunc{E}(\ln |X_n|)$ and $%
\limfunc{Var}(\ln |X_n|)$, let alone to find the distribution of residuals.

\end{document}